\documentclass[12pt]{amsart}
\usepackage{amssymb,latexsym}

\def\ZZ{\mathbb{Z}}
\def\CC{\mathbb{C}}
\def\ee{\mathbf{e}}
\def\PP{\mathbb{P}}
\def\M{\mathcal{M}}

\begin{document}

\title{Quiver Grassmannians and their Euler characteristics: Oberwolfach talk, May 2010}

\author{Andrei Zelevinsky}
\address{\noindent Department of Mathematics, Northeastern University,
  Boston, MA 02115, USA}
\email{andrei@neu.edu}

\begin{abstract}
This is an extended abstract of my talk at the Oberwolfach
Workshop ``Interactions between Algebraic Geometry and Noncommutative Algebra" (May 10 - 14, 2010).
We present some properties of quiver Grassmannians and examples of
explicit computations of their Euler characteristics.
\end{abstract}

\date{June 4, 2010}

\thanks{Research supported in part
by NSF (DMS) grant \# 0801187.}

\maketitle

\makeatletter
\renewcommand{\@evenhead}{\tiny \thepage \hfill  ANDREI ZELEVINSKY \hfill}

\renewcommand{\@oddhead}{\tiny \hfill QUIVER GRASSMANNIANS
 \hfill \thepage}
\makeatother

The aim of this talk is to advertise a very interesting class of algebraic
varieties called \emph{quiver Grassmannians}.
They are defined as follows.
Let $Q$ be a quiver on vertices $\{1, \dots, n\}$.
A \emph{$Q$-representation} is a family $M = (M_i, \varphi_a)$,
where each $M_i$ is a finite dimensional $\CC$-vector space attached to
a vertex~$i$, and each $\varphi_a: M_j \to M_i$ is a linear map attached
to an arrow $a: j \to i$.
The \emph{dimension vector} of~$M$ is the integer vector
${\bf dim}\ M = (\dim M_1, \dots, \dim M_n)$.
A \emph{subrepresentation} of~$M$ is an $n$-tuple of subspaces
$N_i \subseteq M_i$ such that $\varphi_a (N_j) \subseteq N_i$ for
any arrow $a: j \to i$.
With all this terminology in place, for every integer vector
$\ee = (e_1, \dots, e_n)$, the quiver Grassmannian ${\rm
Gr}_\ee(M)$ is defined as the variety of subrepresentations of $M$
with the dimension vector $\ee$.
As a special case, for a one-vertex quiver we get an ordinary
Grassmannian.

Any quiver Grassmannian ${\rm Gr}_\ee(M)$ is Zariski closed
in the product of ordinary Grassmannians $\prod_i {\rm
Gr}_{e_i}(M_i)$, hence is a projective algebraic variety (not
necessarily irreducible or smooth).
Not much is known about their properties.
Motivated by applications to the theory of cluster algebras we
focus on the problem of computing the Euler characteristic $\chi({\rm
Gr}_\ee(M))$.

For a given quiver representation~$M$,
we assemble all the integers $\chi({\rm Gr}_\ee(M))$ into the generating
polynomial $F_M \in \ZZ[u_1, \dots, u_n]$ (\emph{$F$-polynomial}) given by
\begin{equation}
F_M(u_1, \dots, u_n) = \sum_\ee \chi({\rm
Gr}_\ee(M)) u_1^{e_1} \cdots u_n^{e_n} \ .
\end{equation}
It is not hard to show that
$F_{M \oplus N} = F_M F_N$,
hence the study of arbitrary $F$-polynomials $F_M$ reduces to the case of~$M$
indecomposable.

\smallskip

\noindent{\bf Example~1.}
For a one-vertex quiver, i.e., for the ordinary Grassmannians, we have
$F_{\CC^m} = (F_\CC)^m = (1+u)^m$, hence $\chi({\rm
Gr}_e(\CC^m)) = \binom{m}{e}$.

\smallskip

\noindent{\bf Example~2.} Let $Q$ be the \emph{Kronecker quiver}
with two vertices and two arrows from $1$ to $2$.
There are three kinds of indecomposable $Q$-representations:
preprojectives, preinjectives, and regular ones.
More precisely, for every $m \geq 1$, there is a unique (up to an
isomorphism) preprojective indecomposable $Q$-representation $M^{\rm pr}(m)$ of
the dimension vector $(m-1,m)$, a unique preinjective indecomposable
$M^{\rm inj}(m)$ of the dimension vector $(m,m-1)$, and a family of regular
indecomposables $M_\lambda^{\rm reg}(m)$ (parameterized by $\lambda \in \PP^1$)
of the dimension vector $(m,m)$.
As shown in \cite{calzel}, the Euler characteristics of the
corresponding quiver Grassmannians are given as follows:
\begin{align*}
&\chi({\rm Gr}_\ee(M^{\rm pr}(m))) =
\binom{m-e_1}{e_2-e_1}\binom{e_2-1}{e_1},\\
&\chi({\rm Gr}_\ee(M^{\rm inj}(m))) =
\binom{m-e_2}{e_1-e_2}\binom{e_1-1}{e_2}, \\
&\chi({\rm Gr}_\ee(M_\lambda^{\rm reg}(m))) =
\binom{m-e_1}{e_2-e_1}\binom{e_2}{e_1} \ .
\end{align*}

\smallskip

\noindent{\bf Example~3.} Let $Q$ be a Dynkin quiver, i.e., an
orientation of a simply-laced Dynkin diagram; thus, its every connected component
is of one of the $ADE$ types.
The indecomposable $Q$-representations are determined by their
dimension vectors.
Identifying $\ZZ^n$ with the root lattice of the corresponding
root system, these dimension vectors get identified with the
positive roots.
Let $M(\alpha)$ denote the indecomposable representation with the
dimension vector identified with a positive root~$\alpha$.
A unified (type-independent)``determinantal" formula for
$F_{M(\alpha)}$ was given in \cite{yangzel}.
To state it we need some preparation.

First, different orientations of a given Dynkin diagram are in a bijection with
the different Coxeter elements in the Weyl group~$W$:
following \cite{yangzel}, to a Coxeter
element $c = s_{i_1} \cdots s_{i_n}$ (the product of all simple
reflections taken in some order) we associate an orientation with
an edge between $i_\ell$ and $i_k$ oriented from $i_\ell$ to $i_k$
whenever $k < \ell$.

Now let $G$ be the simply-connected semisimple complex algebraic
group associated with our Dynkin diagram.
Each weight $\gamma$ belonging to the $W$-orbit of some
fundamental weight $\omega_i$ gives rise to a regular function
$\Delta_{\gamma, \gamma}$ on $G$ (a \emph{principal generalized
minor}) defined as follows: $\Delta_{\gamma, \gamma}(x)$ is the
diagonal matrix entry of $x \in G$ associated with the one-dimensional
weight subspace of weight~$\gamma$ in the fundamental
representation $V_{\omega_i}$.
Recall also the one-parameter subgroups $x_i(u) = \exp(u e_i)$ and
$y_i(u) = \exp(u f_i)$ in $G$ associated with the Chevalley
generators $e_i$ and $f_i$ of the Lie algebra of~$G$.

With this notation in place we have the following result
essentially proved in \cite{yangzel}:
\begin{equation}
F_{M(\alpha)}(u_1, \dots, u_n) = \Delta_{\gamma, \gamma}
(y_{i_1}(1) \cdots y_{i_n}(1) x_{i_n}(u_{i_n}) \cdots x_{i_1}(u_{i_1})) \ ,
\end{equation}
where $\gamma$ is uniquely determined from the equation
$c^{-1} \gamma - \gamma = \alpha$.

\smallskip

If a quiver $Q$ is \emph{acyclic} (i.e., has no oriented cycles),
the quiver Grassmannians have the following properties:
\begin{itemize}
\item If $M$ is a \emph{general} representation of a given
dimension vector then all its quiver Grassmannians are smooth.
In particular, this is true if $M$ is \emph{rigid}, i.e.,
${\rm Ext}^1(M,M) = 0$ (as explained in \cite[Proposition~3.5]{dwz2},
this follows from the results in \cite{schofield}).
\item If $M$ is indecomposable and rigid then
$\chi({\rm Gr}_\ee(M)) \geq 0$ for all $\ee$.
\end{itemize}

The following example shows that the rigidity condition is essential
for the positivity of the Euler characteristic.

\smallskip

\noindent{\bf Example~4.}
Let $Q$ be the generalized Kronecker quiver with two vertices and
four arrows from $1$ to $2$.
As shown in \cite[Example~3.6]{dwz2}, if $M$ is a general
representation of dimension vector $(3,4)$, and $\ee = (1,3)$
then ${\rm Gr}_\ee(M)$ is isomorphic to a smooth projective curve
of degree~$4$ in $\PP^2$, hence $\chi({\rm Gr}_\ee(M)) = -4$.

\smallskip

With applications to the theory of
cluster algebras in mind, we will assume from now on
that every quiver $Q$ under consideration has no
loops or oriented $2$-cycles.
If in addition a quiver $Q$ is acyclic, the quiver Grassmannians most important for
these applications are those in rigid indecomposable $Q$-representations.
In \cite{dwz2} it was suggested how to extend this class of
representations to the case of not necessarily acyclic quivers.
Here is a brief account of this approach.

First, we restrict our attention to $Q$-representations satisfying
relations imposed by a generic \emph{potential} $S$ on $Q$.
Roughly speaking, $S$ is a generic linear combination (possibly
infinite) of cyclic paths in $Q$, viewed as an element of the
completed path algebra $\widehat{\CC Q}$ (see \cite{dwz} for a detailed setup).
By a $(Q,S)$-representation we mean a
$Q$-representation annihilated by sufficiently long
paths in $\CC Q$ and by all \emph{cyclic derivatives} of~$S$.

Second, following \cite{mrz,dwz}, we consider \emph{decorated} $(Q,S)$-representations:
such a representation is a pair $\M = (M,V)$,
where $M$ is a $(Q,S)$-representation as above, and $V = (V_i)$
is a collection of finite-dimensional vector spaces attached
to the vertices (with no maps attached).

\emph{Mutations} of quivers with potentials and their
representations at arbitrary vertices were introduced and studied in \cite{dwz}.
As shown in \cite{dwz2}, a natural class of
$(Q,S)$-representations generalizing rigid indecomposable
representations of acyclic quivers consists of
$(Q,S)$-representations obtained by mutations from \emph{negative
simple} representations (those having some $V_i$ equal to $\CC$, and the rest
of the spaces $M_j$ and $V_j$ equal to~$0$).
One of the main constructions in \cite{dwz2} is that of an
integer-valued function $E(\M)$ on the set of
$(Q,S)$-representations, which is invariant under mutations and
vanishes on negative simple representations.
An open question is whether the condition $E(\M) = 0$ for an
indecomposable $(Q,S)$-representation $\M$ implies that $\M$ can be
obtained by mutations from a negative simple representation.
Regardless of the answer to this question, quiver Grassmannians in
the indecomposable $(Q,S)$-representations with $E(\M) = 0$ deserve further study.

\end{document}